\documentclass[english]{amsart}
\usepackage[T1]{fontenc}
\usepackage[latin9]{inputenc}
\usepackage{amssymb}

\makeatletter

\providecommand{\tabularnewline}{\\}

\numberwithin{equation}{section} 
\numberwithin{figure}{section} 

\newtheorem{Theorem}{Theorem}[section]
\newtheorem{Cor}[Theorem]{Corollary}
\newtheorem{Lemma}[Theorem]{Lemma}
\newtheorem{Proposition}[Theorem]{Proposition}
\theoremstyle{definition}
\newtheorem{Definition}[Theorem]{Definition}
\theoremstyle{remark}

\numberwithin{equation}{section}

\newcommand{\R}{\mathbb R}
\newcommand{\N}{\mathbb N}
\newcommand{\Z}{\mathbb Z}

\newcommand{\C}{\mathcal{C}}
\newcommand{\B}{\mathcal{B}}
\newcommand{\p}{\mathcal{P}}
\newcommand{\A}{\mathcal{A}}
\newcommand{\rel}{\mathcal{R}}

\newcommand{\F}{\mathcal{F}}

\newcommand{\lo}{\mathcal{L}}
\newcommand{\Hol}{\mathcal{H}}
\def\d{\partial}
\def\K{{ \! \rm \ I\!K}}

\makeatother

\usepackage{babel}

\begin{document}

\title{Ambrose-Singer theorem on diffeological bundles and complete integrability of the KP equation}

\author{Jean-Pierre Magnot}

\address{Lyc\'ee Jeanne d'Arc - 40 avenue de Grande Bretagne - F63000 Clermont-Ferrand}

\email{jean-pierr.magnot@ac-clermont.fr}

\begin{abstract}In this paper, we start from an extension of 
the notion of holonomy on diffeological bundles, reformulate 
the notion of regular Lie group or Fr\"olicher Lie groups, state 
an Ambrose-Singer theorem that enlarges the one stated in \cite{Ma2}, and 
conclude with a differential geometric treatment of KP hierarchy. 
The examples of Lie groups that are studied 
are principally those obtained by enlarging some graded Fr\"olicher (Lie) algebras 
such as formal $q-$series of the quantum algebra of pseudo-differential operators. These deformations can be defined for classical 
pseudo-differential operators but they are used here on 
formal pseudo-differential operators in order to get a differential
geometric framework to deal with the KP hierarchy that is known
to be completely integrable in the sense of Frobenius. 
Here, we get an integration of the Zakharov-Shabat connection form
by means of smooth sections of a (differential geometric) bundle with structure 
group some groups of $q-$deformed operators. The integration obtained 
by Mulase \cite{M1}, and the key tools he developped, are totally recovered on the germs of the smooth maps of our construction. 

The tool coming from (classical) differential geometry used in this 
construction is the holonomy group, on which we have an 
Ambrose-Singer-like theorem: the Lie algebra is spanned 
by the curvature elements. This result is proved for any connection 
a diffeological principal bundle with structure group a regular Fr\"olicher Lie group. The case of a (classical) Lie group modelled on a complete locally convex topological vector space is also recovered and the work developped in \cite{Ma,Ma2} is completed.

\end{abstract}
\maketitle


\noindent
{\small MSC (2010) : 58B25; 58Z05; 37K10; 37K25; 37K30}

\noindent
{\small Keywords : KP hierarchy; holonomy; Ambrose-Singer theorem; Infinite dimensional Lie groups; Fr\"olicher Lie groups; Diffeological spaces}

\tableofcontents
\section*{Introduction}
The aim of this paper is to give some new results of integration of horizontal lifts of paths on a diffeological bundle that can be applied to integration of systems or hierarchies of equations that cannot be expressed in the standard framework of Lie groups, such as the KP hierarchy. 

For this, we have to develop the same tools on a wider category: the category of diffeological spaces. This category, described first by Souriau (see e.g. \cite{Sou}), is very easy to use on infinite catersian products but carries intrinsic theoretical ambiguities. For example, even for compact manifolds, there exists many diffeologies underlying the same manifold structure, see e.g. section 1.4. This is why we felt the need in \cite{Ma} 
to compare it to the notion of 
Fr\"olicher space \cite{FK}, which appears 
to complete the first one (this idea has been more recently 
developped in \cite{Wa}, in a more theoretical approach). 
This framework seems
to us very useful to study 
infinite systems of equations, and 
we show that it enables one to use 
fully justified differential geometric tools ( here,  an 
Ambrose-Singer theorem) 
in the context of the Kadomtsev-Petviashvili hierarchy
already discussed by Mulase \cite{M1,M3} in a more algebraic way.
   
\vskip 12pt
For this,

- We first end the study begun in \cite{Ma,Ma2} 
on the Ambrose-Singer like theorems that 
one can state out of the settings of 
finite dimensional manifolds.  
We began in \cite{Ma2} by stating an Ambrose-Singer theorem for a class of infinite dimensional 
principal bundles, addressing also open questions on the mathematical structure of the 
holonomy groups constructed. In this work, completed in \cite{Ma}, 
the most natural frameworks are not manifolds but weaker frameworks such as diffeological or Fr\"olicher spaces. 
The starting point of \cite{Ma2} was a Lie theorem stated by Robart \cite{Rob} for
a wide class of Lie algebras, which helped us to circumvent the lack of easy-to-use 
Frobenius theorem on integrable distributions. 
The approach is here slightly changed. The objects considered here are as general as possible: 
diffeological or Fr\"olicher spaces. Infinite dimensional manifolds appear as particular cases.
The notion of connection is constructed step by step, 
using an ``integrated'' notion of connections on diffeological bundles from the
``infinitesimal'' approach of Iglesias \cite{Ig}. 
This approach is adapted here to define path-liftings, a generalization 
of horizontal lifts of paths in section \ref{3} which enables us to define holonomy groups. 
Connection forms are the infinitesimal aspect of $G-$invariant path-liftings. This leads us to 
a classification of diffeological structure of a Fr\"olicher space which is less subtle that the one described in \cite{Wa} and initiated by Iglesias \cite{Igdiff} (but sufficient for our purpose), 
a Lie theorem for Fr\"olicher Lie subalgebras (Theorem \ref{Lie}) and the definition of 
regular Fr\"olicher Lie groups with regular Lie algebra 
which is directly derived from on the approach of Leslie 
\cite{Les} for analogous notions on diffeological groups. 
These notions are applied to very easy examples: 
generalized Lie groups defined by Omori in \cite{Om} 
(Proposition \ref{omo}), and "graded" Fr\"olicher algebras $\A$ 
that appear to be enlargeable into regular Fr\"olicher 
Lie groups (Theorem \ref{regulardeformation}).
All these results are technical tools 
which allows the machinery of \cite{Ma2} to work in the desired settings: we can build a second holonomy group, 
noted $H^{red},$ which is the smaller in a certain class to which we can apply the theorem of 
reduction of the structure group \ref{Courbure}.
We have in this setting the following Ambrose-Singer theorem: the Lie algebra of 
$H^{red}$ is the regular Lie algebra spanned by the curvature elements, 
that is the smaller regular Lie algebra which contains the curvature elements.

- Secondly, an easy application of the Ambrose-Singer theorem is the construction 
of smooth sections of principal bundles with contractible basis via 
flat connections. We apply this procedure to deal with the complete integrability of the Kadomtsev-Petviashvili equation (or KP equation for short)
which is known to be an equation equivalent to a $0-$currvature condition 
on the so-called Zhakharov-Shabat connection 
(See e.g. \cite{M1,Ue} for a standard description). 
The variables the hierarchy are $(t_1,t_2,...)$ so that the base of this principal bundle is the (algebraic) $k-$vector space $\oplus_{n \in \N^*} k t_n,$
which can be viewed as a Fr\"olicher space. 
The integration of this equation 
has been (algebraically) given by Mulase \cite{M1} in the context of 
formal power series with an infinite number of variables, using various results given by other authors and a so-called Birkhoff decomposition on groups of formal power series.
These groups of formal series are understood as formal expressions for global sections of a principal bundle. We show here that the KP hierarchy is equivalent 
to another hierarchy (one should say $q-$deformed hierarchy, but we have to say that its seems to have no straightforward correspondence with a problem related to deformation quantization), which has an analogous 0-curvature condition. 
The structure groups of the bundles considered are regular. These are groups given as examples of regular Fr\"olicher Lie groups in section \ref{1.}, namely groups obtained from algebras of formal pseudo-differential operators. This enables us to fully apply the Ambrose-Singer theorem announced before, and recover a Birkoff-type decomposition by this way. This decomposition fits with the one 
described by Mulase by taking the germs of the smooth sections we obtained, and setting the deformation parameter $q=1.$

\vskip 12pt
\noindent
As we said at the beginning of the introduction, 
we think that such an approach can be very promising 
for other equations or hierarchies 
based on formal series. The first example that seem similar is the super 
KP equation \cite{M2} and we hope to be able very soon 
to give a result on this system of equation.

\section{Preliminaries on differentiable structures}

\label{1.} The objects of the category of -finite or infinite- dimensional
smooth manifolds is made of topological spaces $M$ equipped with
a collection of charts called maximal atlas that enables one to make
differentiable calculus. But for 
e.g. projective limits of manifolds or infinite products,
a differential calculus is needed where as no atlas can be defined.
To circumvent this problem which occurs in various frameworks, several
authors have independently developed some ways to define differentiation
without defining charts. We use here three of them. The first one
is due to Souriau \cite{Sou}, the second one is due to Sikorski,
and the third one is a setting closer to the setting of differentiable
manifolds is due to Frölicher (see e.g. \cite{CN} for an introduction
on these two last notions). In this section, we review some basics
on these three notions.

\subsection{Souriau's diffeological spaces, Sikorski's differentiable spaces,
Frölicher spaces}
\label{1.1}

\begin{Definition} Let $X$ be a set.

\noindent $\bullet$ A \textbf{plot} of dimension $p$ (or $p$-plot)
on $X$ is a map from an open subset $O$ of $\R^{p}$ to $X$.

\noindent $\bullet$ A \textbf{diffeology} on $X$ is a set $\p$
of plots on $X$ such that, for all $p\in\N$,

- any constant map $\R^{p}\rightarrow X$ is in $\p$;

- Let $I$ be an arbitrary set; let $\{f_{i}:O_{i}\rightarrow X\}_{i\in I}$
be a family of maps that extend to a map $f:\bigcup_{i\in I}O_{i}\rightarrow X$.
If $\{f_{i}:O_{i}\rightarrow X\}_{i\in I}\subset\p$, then $f\in\p$.

- (chain rule) Let $f\in\p$, defined on $O\subset\R^{p}$. Let $q\in\N$,
$O'$ an open subset of $\R^{q}$ and $g$ a smooth map (in the usual
sense) from $O'$ to $O$. Then, $f\circ g\in\p$.

\vskip 6pt $\bullet$ If $\p$ is a diffeology $X$, $(X,\p)$ is
called \textbf{diffeological space}.

\noindent Let $(X,\p)$ et $(X',\p')$ be two diffeological spaces,
a map $f:X\rightarrow X'$ is \textbf{differentiable} (=smooth) if
and only if $f\circ\p\subset\p'$. \end{Definition}

\noindent \textbf{Remark.} Notice that any diffeological space $(X,\p)$
can be endowed with the weaker topology such that all the maps that
are in $\p$ are continuous. But we prefer to mention this only for
memory as well as other questions that are not closely related to
our construction, and stay closer to the goals of this paper.

Let us now define the Sikorski's differential spaces. Let $X$ be
a Haussdorf topological space.

\begin{Definition}

\noindent $\bullet$ A (Sikorski's) \textbf{differential space} is
a pair $(X,\F)$ where $\F$ is a family of maps $X\rightarrow\R$
such that

- the topology of $X$ is the initial topology with respect to $\F$

- for any $n\in\N$, for any smooth map $\varphi:\R^{n}\rightarrow\R$,
for any $(f_{1},...,f_{n})\in\F^{n}$, $\varphi\circ(f_{1},...,f_{n})\in\F$.

\noindent $\bullet$ Let $(X,\F)$ et $(X',\F')$ be two differential
spaces, a map $f:X\rightarrow X'$ is \textbf{differentiable} (=smooth)
if and only if $\F'\circ f\subset\F.$ \end{Definition}

We now introduce Frölicher spaces.

\begin{Definition} $\bullet$ A \textbf{Frölicher} space is a triple
$(X,\F,\C)$ such that

- $\C$ is a set of paths $\R\rightarrow X$,

- A function $f:X\rightarrow\R$ is in $\F$ if and only if for any
$c\in\C$, $f\circ c\in C^{\infty}(\R,\R)$;

- A path $c:\R\rightarrow X$ is in $\C$ (i.e. is a \textbf{contour})
if and only if for any $f\in\F$, $f\circ c\in C^{\infty}(\R,\R)$.

\vskip 5pt $\bullet$ Let $(X,\F,\C)$ et $(X',\F',\C')$ be two
Frölicher spaces, a map $f:X\rightarrow X'$ is \textbf{differentiable}
(=smooth) if and only if $\F'\circ f\circ\C\in C^{\infty}(\R,\R)$.
\end{Definition}

Any family of maps $\F_{g}$ from $X$ to $\R$ generate a Frölicher
structure $(X,\F,\C)$, setting \cite{KM}:

- $\C=\{c:\R\rightarrow X\hbox{ such that }\F_{g}\circ c\subset C^{\infty}(\R,\R)\}$

- $\F=\{f:X\rightarrow\R\hbox{ such that }f\circ\C\subset C^{\infty}(\R,\R)\}.$

One easily see that $\F_{g}\subset\F$. This notion will be useful
in the sequel to describe in a simple way a Frölicher structure.

A Frölicher space, as a differential space, carries a natural topology,
which is the pull-back topology of $\R$ via $\F$. In the case of
a finite dimensional differentiable manifold, the underlying topology
of the Frölicher structure is the same as the manifold topology. In
the infinite dimensional case, these two topologies differ very often.

\vskip 12pt

In the three previous settings, we call $X$ a \textbf{differentiable
space}, omitting the structure considered. Notice that, in the three
previous settings, the sets of differentiable maps between two differentiable
spaces satisfy the chain rule. Let us now compare these three settings:
One can see (see e.g. \cite{CN}) that we have the following, given
at each step by forgetful functors: \vskip 10pt \begin{tabular}{ccccc}
smooth manifold  & $\Rightarrow$  & Frölicher space  & $\Rightarrow$  & Sikorski differential space\tabularnewline
\end{tabular}\vskip 10pt Moreover, one remarks easily from the definitions that,
if $f$ is a map from a Frölicher space $X$ to a Frölicher space
$X'$, $f$ is smooth in the sense of Frölicher if and only if it
is smooth in the sense of Sikorski.

\vskip 12pt

One can remark, if $X$ is a Frölicher space, we define a natural
diffeology on $X$ by \cite{Ma}: 
$$
\p(\F)=
\coprod_{p\in\N}\{\, f\hbox{ p-
paramatrization on } X; \, \F \circ f \in C^\infty(O,\R) \quad \hbox{(in
the usual sense)}\}.$$
With this construction, we also get a natural diffeology when
$X$ is a Frölicher space. In this case, one can easily show the following:
\begin{Proposition}\cite{Ma}
Let$(X,\F,\C)$
and $(X',\F',\C')$ be two Frölicher spaces. A map $f:X\rightarrow X'$
is smooth in the sense of Frölicher if and only if it is smooth for
the underlying diffeologies. \end{Proposition}

Thus, we can also state:

\begin{tabular}{ccccc}
smooth manifold  & $\Rightarrow$  & Frölicher space  & $\Rightarrow$  & Diffeological space\tabularnewline
\end{tabular}

For a more complete description of the relationship between
diffeological spaces and Fr\"olicher spaces, see \cite{Wa}.

\subsection{Push-forward, quotient and trace}

We give here only the results that will be used in the sequel.

\begin{Proposition} \cite{Ma} Let $(X,\p)$ be a diffeological space,
and let $X'$ be a set. Let $f:X\rightarrow X'$ be a surjective map.
Then, the set \[
f(\p)=\{u\hbox{such that}u\hbox{restricts to some maps of the type}f\circ p;p\in\p\}\]
 is a diffeology on $X'$, called the \textbf{push-forward diffeology}
on $X'$ by $f$. \end{Proposition}

We have now the tools needed to describe the diffeology on a quotient:

\begin{Proposition} \label{quotient} let $(X,\p)$ b a diffeological
space and $\rel$ an equivalence relation on $X$. Then, there is
a natural diffeology on $X/\rel$, noted by $\p/\rel$, defined as
the push-forward diffeology on $X/\rel$ by the quotient projection
$X\rightarrow X/\rel$. \end{Proposition}

Given a subset $X_{0}\subset X$, where $X$ is a Frölicher space
or a diffeological space, we can define on trace structure on $X_{0}$,
induced by $X$.

$\bullet$ If $X$ is equipped with a diffeology $\p$, we can define
a diffeology $\p_{0}$ on $X_{0}$ setting \[
\p_{0}=\lbrace p\in\p\hbox{such that the image of }p\hbox{ is a subset of }X_{0}\rbrace.\]

$\bullet$ If $(X,\F,\C)$ is a Frölicher space, we take as a generating
set of maps $\F_{g}$ on $X_{0}$ the restrictions of the maps $f\in\F$.
In that case, the contours (resp. the induced diffeology) on $X_{0}$
are the contours (resp. the plots) on $X$ which image is a subset
of $X_{0}$.

\subsection{Cartesian products and projective limits}

The category of Sikorski differential spaces is not cartesianly closed,
see e.g. \cite{CN}. This is why we prefer to avoid the questions
related to cartesian products on differential spaces in this text,
and focuse on Frölicher and diffeological spaces, since the cartesian
product is a tool essential for the definition of configuration spaces.

In the case of diffeological spaces, we have the following \cite{Sou}:

\begin{Proposition} \label{prod1} Let $(X,\p)$ and $(X',\p')$
be two diffeological spaces. We call \textbf{product diffeology} on
$X\times X'$ the diffeology $\p\times\p'$ made of plots $g:O\rightarrow X\times X'$
that decompose as $g=f\times f'$, where $f:O\rightarrow X\in\p$
and $f':O\rightarrow X'\in\p'$. \end{Proposition}

Then, in the case of a Frölicher space, we derive very easily, compare
with e.g. \cite{KM}: \begin{Proposition} \label{prod2} Let $(X,\F,\C)$
and $(X',\F',\C')$ be two Frölicher spaces, with natural diffeologies
$\p$ and $\p'$ . There is a natural structure of Frölicher space
on $X\times X'$ which contours $\C\times\C'$ are the 1-plots of
$\p\times\p'$. \end{Proposition}

We can even state the same results in the case of infinite products,
in a very trivial way by taking the cartesian products of the plots or of the contours.  Let us now give the description of what happens
for projective limits of Frölicher and diffeological spaces.

\begin{Proposition} \label{froproj} Let $\Lambda$ be an infinite
set of indexes, that can be uncoutable.

$\bullet$ Let $\lbrace(X_{\alpha},\p_{\alpha})\rbrace_{\alpha\in\Lambda}$
be a family of diffeological spaces indexed by $\Lambda$ totally
ordered for inclusion,  with $(i_{\beta,\alpha})_{(\alpha,\beta)\in\Lambda^{2}}$
a family of diffeological maps . If $X=\bigcap_{\alpha\in\Lambda}X_{\alpha},$
$X$ carries the \textbf{projective diffeology} $\p$ which is the
pull-back of the diffeologies $\p_{\alpha}$ of each $X_{\alpha}$
via the family of maps $(f_{\alpha})_{\alpha\in\Lambda}.$ The diffeology
$\p$ made of plots $g:O\rightarrow X$ such that, for each $\alpha\in\Lambda,$
\[
f_{\alpha}\circ g\in\p_{\alpha}.\]
 This is the biggest diffeology for which the maps $f_{\alpha}$ are
smooth.

$\bullet$ Let $\lbrace(X_{\alpha},\F_{\alpha},\C_{\alpha})\rbrace_{\alpha\in\Lambda}$
be a family of Frölicher spaces indexed by $\Lambda$ totally ordered
for inclusion,  with $(i_{\beta,\alpha})_{(\alpha,\beta)\in\Lambda^{2}}$
a family of differentiable maps . with natural diffeologies $\p_{\alpha}.$
There is a natural structure of Frölicher space $X=\bigcap_{\alpha\in\Lambda}X_{\alpha},$
which contours \[
\C=\bigcap_{\alpha\in\Lambda}\C_{\alpha}\]
 are the 1-plots of $\p=\bigcap_{\alpha\in\Lambda}\p_{\alpha}.$ A
generating set of functions for this Frölicher space is the set of
maps of the type: 
$$
\bigcup_{\alpha \in \Lambda} \F_{\alpha} \circ f_{\alpha}.$$

\end{Proposition}

\subsection{Differential forms on a diffeological space and differential dimension.}

\begin{Definition} \cite{Sou}
Let $(X,\p)$ be a diffeological space and let $V$ be a vector space equipped with a differentiable structure. 
A $V-$valued $n-$differential form $\alpha$ on $X$ (noted $\alpha \in \Omega^n(X,V))$ is a map 
$$ \alpha : \{p:O_p\rightarrow X\} \in \p \mapsto \alpha_p \in \Omega^n(p;V)$$
such that 

$\bullet$ Let $x\in X.$ $\forall p,p'\in \p$ such that $x\in Im(p)\cap Im(p')$, 
the forms $\alpha_p$ and $\alpha_{p'}$ are of the same order $n.$ 

$\bullet$ Moreover, let $y\in O_p$ and $y'\in O_{p'}.$ If $(X_1,...,X_n)$ are $n$ germs of paths in 
$Im(p)\cap Im(p'),$ if there exists two systems of $n-$vectors $(Y_1,...,Y_n)\in (T_yO_p)^n$ and $(Y'_1,...,Y'_n)\in (T_{y'}O_{p'})^n,$ if $p_*(Y_1,...,Y_n)=p'_*(Y'_1,...,Y'_n)=(X_1,...,X_n),$
$$ \alpha_p(Y_1,...,Y_n) = \alpha_{p'}(Y'_1,...,Y'_n).$$

We note by $$\Omega(X;V)=\oplus_{n\in \mathbb{N}} \Omega^n(X,V)$$ the set of $V-$valued differential forms.  
\end{Definition}

With such a definition, we feel the need to make two remarks for the reader:
 
$\bullet$ If there does not exist $n$ linearly independent vectors $(Y_1,...,Y_n)$
defined as in the last point of the definition, $\alpha_p = 0$ at $y.$

$\bullet$ Let $(\alpha, p, p') \in \Omega(X,V)\times \p^2.$ 
If there exists $g \in C^\infty(D(p); D(p'))$ (in the usual sense) 
such that $p' \circ g = p,$ then $\alpha_p = g^*\alpha_{p'}.$ 

\vskip 12pt
\begin{Proposition}
The set $\p(\Omega^n(X,V))$ made of maps $q:x \mapsto \alpha(x)$ from an open subset $O_q$ of a 
finite dimensional vector space to $\Omega^n(X,V)$ such that for each $p \in \p,$ $$\{ x \mapsto \alpha_p(x) \} \in C^\infty(O_q, \Omega^n(O_p,V)),$$
is a diffeology on $\Omega^n(X,V).$  
\end{Proposition}

Working on plots of the diffeology, one can define the product and the differential of differential forms, which have the same properties as  
the product and the differential of differential forms. 

\begin{Definition}
Let $(X,\p)$ be a diffeological space. 

\noindent
$\bullet$ $(X,\p)$ is \textbf{finite-dimensional} if and only if $$\exists n_0\in\mathbb{N},\quad \forall n\in \mathbb{N}, \quad n\geq n_0 \Rightarrow dim(\Omega^n(X,\mathbb{R}))=0.$$
Then, we set $$dim(X,\p)=max\{n\in \mathbb{N}| dim(\Omega^n(X,\mathbb{R}))>0\}.$$
\noindent
$\bullet$ If not, $(X,\p)$ is called \textbf{infinite dimensional}.
\end{Definition}

Let us make a few remarks on this definition.
If $X$ is a manifold with $dim(X)=n,$ the natural diffeology as described in section \ref{1.1} 
(also called ``n\'ebuleuse'' diffeology) is such that $$dim(X,\p_0)=n.$$
Now, if $(X,\F,\C)$ is the natural Fr\"olicher structure on $X,$ take $\p_1$ generated by the maps of the type
$g\circ c$, where $c\in \C$ and $g$ is a smooth map from an open subset of a finite dimensional space to $\mathbb{R}.$
This is an easy exercise to show that $$dim(X,\p_1)=1.$$
This first point shows that the dimension depends on the diffeology considered. Now, we remark that $\F$
is the set of smooth maps $(X,\p_1)\rightarrow \mathbb{R},$

This leads to the following definition, since $\p(\F)$ is clearly the diffeology with the biggest dimension associated to $(X,\F,\C)$:

\begin{Definition}
The \textbf{dimension} of a Fr\"olicher space $(X,\F,\C)$ is the dimension of the diffeologial space $(X,\p(\F)).$
\end{Definition}

\subsection{Regular Frölicher groups}

Let $(G, \F, \C)$ be a Fr\"'olicher space which is a group  such that
the group law and the inversion map are smooth. These laws are also smooth
for the underlying diffeology. Then, following \cite{Les}, this is
possible as in the case of manifolds to define a tangent space and
a Lie algebra $\mathfrak{g}$ of $G$ using germs of smooth maps.
Let us precise the algebraic, diffeological and Fr\"olicher structures of $\mathfrak{g}.$

\begin{Proposition}
Let $\mathfrak{g} = \{ \partial_t c(0) ; c \in \C \hbox{ and } c(0)=e_G \}$
be the space of germs of paths at $e_G.$
\begin{itemize}
	\item Let $(X,Y) \in \mathfrak{g}^2,$ $X+Y = \partial_t(c.d)(0)$  where $c,d \in \C ^2,$ $c(0) = d(0) =e_G ,$ 
	$X = \partial_t c(0)$ and $Y = \partial_t d(0).$ 
	\item Let $(X,g) \in \mathfrak{g}\times G,$ $Ad_g(X) = \partial_t(g c g^{-1})(0)$  where $c \in \C ,$ $c(0) =e_G ,$ 
	and $X = \partial_t c(0).$  
	\item Let $(X,Y) \in \mathfrak{g}\times G,$ $[X,Y] = \partial_t( Ad_{c(t)}Y)$   where $c \in \C ,$ $c(0) =e_G ,$ 
	$X = \partial_t c(0).$
\end{itemize}
All these operations are smooth and thus well-defined. 
\end{Proposition} 

The basic properties remain globally the same as in the case of Lie groups, and the prrofs are similar replacing charts by plots of the underlying diffeologies. Since the work is already done in \cite{Les} in the category of diffeological groups, we refer the reader to it for details. 

\begin{Definition} A Frölicher group $G$ with Lie algebra $\mathfrak{g}$
is called \textbf{regular} if and only if there is a smooth map \[
Exp:C^{\infty}([0;1],\mathfrak{g})\rightarrow C^{\infty}([0,1],G)\]
 such that $g(t)=Exp(v(t))$ if and only if $g$ is the unique solution
of the differential equation \[
\left\{ \begin{array}{l}
g(0)=e\\
\frac{dg(t)}{dt}g(t)^{-1}=v(t)\end{array}\right.\]
 We define \begin{eqnarray*}
exp:\mathfrak{g} & \rightarrow & G\\
v & \mapsto & exp(v)=g(1)\end{eqnarray*}
 where $g$ is the image by $Exp$ of the constant path $v.$ \end{Definition}

The classical setting for infinite dimensional differential geometry requires the model topological vector space
to be complete or Mac-Key complete. One of the reasons for this choice 
is to ensure the existence of the integral of a path over a 
compact interval. This means that the choice of an adaquate topology is 
necessary. For vector spaces, such a study can be found 
in \cite{KM}, when the properties of the so-called ``convenient vector spaces''
are given. We have to remark that a vector space for which addition and scalar multiplication
are compatible with a given Fr\"olicher structure needs only 
a topological structure to become a convenient vector space.
In order to circumvent these topological considerations, and adapting the     
terminology of regular Lie groups to vector spaces (which are viewed as abelian Lie groups), we set:

\begin{Definition}
Let $(V,\F, \C)$ be a Fr\"olicher vector space, i.e. a vector space $V$ equipped with a Fr\"'olicher structure compatible
with the vector space addition and the scalar multiplication. $(V,\F, \C)$ is \textbf{regular} if there is a smooth map 
$$ \int_0^{(.)} : C^\infty([0;1];V) \rightarrow C^\infty([0;1],V)$$ such that $\int_0^{(.)}v = u$ if and only if $u$ is the unique solution of 
the differential equation
\[
\left\{ \begin{array}{l}
u(0)=0\\
u'(t)=v(t)\end{array}\right. .\]

\end{Definition}

This definition is of course fulfilled if $V$ is a complete locally convex topological vector space, equipped with its natural Fr\"olicher structure.

\begin{Definition}
Let $G$ be a Fr\"olicher Lie group with Lie algebra $\mathfrak{g}.$ Then, $G$ is \textbf{regular with regular Lie algebra}
if both $G$ and $\mathfrak{g}$ are regular.
\end{Definition}
Following the terminology used in the early inverstigations on infinite dimensional Lie theory, we say that a regular Lie algebra $\mathfrak{g}$ is \textbf{enlargeable}
if there exists a Fr\"olicher Lie group $ G$ 
with Lie algebra   $\mathfrak{g}.$ 
\begin{Theorem}
\label{Lie}
Let $G$ be a regular Fr\"olicher Lie group with Lie algebra $\mathfrak{g}.$
Let $\mathfrak{g}_1$ be a Lie subalgebra of $\mathfrak{g}.$
Let $G_1 = Exp(C^\infty([0;1];\mathfrak{g}_1))(1).$
If $Ad_{G_1 \cup G_1^{-1}}(\mathfrak{g_1})=\mathfrak{ g}_1,$ $G_1$ is a Fr\"olicher subgroup of $G.$
\end{Theorem}

\noindent
\textbf{Proof.}
$G_1$ is obviously a Fr\"olicher subspace of $G.$ 
All we need to show is that $G_1$ is a subgroup of $G$ (algebraically).
This is a well-known procedure, found by Robart \cite{Rob} to our knowledge:

 - If $g= Exp(u(.))(1)\in G_1$ and $g' = Exp(v(.))(1)\in G_1,$ with $u$ and $v$ 
smooth paths that
are stationary at the endpoints, $gg' = Exp(Ad_g(v)\vee u)\in G_1.$

- If $g= Exp(u(.))(1)$ and $g(t) = Exp(u(.))(t),$ $g^{-1}=Exp(-Ad_{g^{-1}}u(.))(1) \in G_1.$ \qed

\vskip 12pt
We call such a subalgebra $\mathfrak{g}_1$ an $Ad-$\textbf{stable} Lie subalgebra. 
\begin{Lemma}
Let $\mathfrak{g}_0 $ be a Lie subalgebra of an enlargeable regular Fr\"olicher Lie algebra $\mathfrak{g}.$ Then there exists a $Ad-$stable regular Fr\"olicher Lie subalgebra $\mathfrak{g}_1 \subset \mathfrak{g},$ minimal for inclusion. 
\end{Lemma}

\vskip 12pt
\noindent
\textbf{Proof.}
Let $E$ be the set of $Ad-$stable regular Fr\"olicher Lie subalgebras of 
$\mathfrak{g} $ that contain $\mathfrak{g}_0.$  We have $\mathfrak{g} \in E$ so that $E$ is nonempty. 
Let $I$ be an index est and $\{\mathfrak{g}_i\}_{i \in I} \in E^I.$ Let $G_i=Exp(C^\infty([0;1];\mathfrak{g}_i))(1).$
The condition $Ad_{G_i \cup G_i^{-1}}(\mathfrak{g}_i)=\mathfrak{ g}_i$
is equivalent to 
$Ad_{G_i \cup G_i^{-1}}(\mathfrak{g}_i)\subset \mathfrak{ g}_i.$
So that, if $\{\mathfrak{g}_i\} \in E^I,$ if  $G'=Exp(C^\infty([0;1];\cap_{i \in I}\mathfrak{g}_i))(1),$
we have that \begin{eqnarray*}Ad_{G' \cup G'^{-1}}(\cap_{i \in I} \mathfrak{g}_i) &\subset& \bigcap_{i \in I} Ad_{G' \cup G'^{-1}}(\mathfrak{g}_i)\\
& \subset &  \bigcap_{i \in I} Ad_{G_i \cup G_i^{-1}}(\mathfrak{g}_i)  \\
& \subset & \bigcap_{i \in I}\mathfrak{g}_i.
\end{eqnarray*}
Since $E$ is stable for intersection and since $\mathfrak{g}_0$ is contained in any elent of $E,$ the set $E$ has an element $\mathfrak{g}_1\neq \emptyset,$ minimal for inclusion. \qed

\begin{Definition}
with the previous notations, we call $\mathfrak{g}_1$ the regular Lie algebra $Ad-$\textbf{spanned} by $\mathfrak{g}_0.$
\end{Definition}

We now turn to the examples, and we feel the need to begin with Omori's generalized Lie groups. 

\begin{Proposition} \label{omo} Let $(G_{n})_{n\in\N}$ be a sequence of Banach
Lie groups, increasing for $\supset,$ and such that the inclusions
are Lie group morphisms. Let $G=\bigcap_{n\in\N}G_{n}.$ Then, $G$
is a Frölicher regular Lie group with regular Lie algebra $\mathfrak{g}=\bigcap_{n\in\N}\mathfrak{g}_{n}.$
\end{Proposition} 

\textbf{Proof.} In this proof, each Lie group $G_{n}$
is equipped with its natural Frölicher structure of smooth manifold
$(G_{n},\F_{n},\C_{n}),$ with underlying diffeology $\p_{n}.$ The
group $G$ is equipped with the projective Frölicher structure $(G,\F,\C),$
with underlying diffeology $\p.$ Let $f,g\in\p^{2}$ and $\alpha\in\Lambda.$
Then, $f_{\alpha}\circ(f.g)=(f_{\alpha}\circ f).(f_{\alpha}\circ g)\in\p_{\alpha}$
and $f_{\alpha}\circ g^{-1}=(f_{\alpha}\circ g)^{-1}\in\p_{\alpha}.$
Thus, multiplication and inversion are differentiable in $G,$ in
the sense of diffeologies and hence in the sense of Frölicher. $G$
is a Frölicher Lie group.

\noindent We now look at $\mathfrak{g},$ which equals to $\bigcap_{n\in\N}\mathfrak{g}_{n}$
trivially. Then, for each $n\in\N,$ if $Exp_{n}$ is the exponential
on $G_{n}$, $\forall m<n,Exp_{m}\circ D_{e}i_{m,n}=Exp_{n}$ since
$i_{m,n}$ is a morphism of Lie groups. (in fact, if we want to be
rigorous, we need to replace $D_{e}i_{m,n}$ by the map $C^{\infty}([0;1];\mathfrak{g}_{n})\rightarrow C^{\infty}([0;1];\mathfrak{g}_{m})$
induced by $i_{m,n}$.) Thus,the exponential on $G$ is smooth. \qed
\vskip 12pt

Let us now give another class of examples of regular Lie groups. 

\begin{Theorem} \label{regulardeformation}
Let $(A_n)_{n \in \N^*} $ be a sequence of complete locally convex (Fr\"olicher)
vector spaces which are regular, 
equipped with a graded smooth multiplication operation
on $ \bigoplus_{n \in \N^*} A_n ,$ i.e. a multiplication such that 
$A_n .A_m \subset A_{n+m},$ smooth with respect to the corresponding Fr\"olicher structures.
Let us assume that:

Then, the set 
$$1 + \A = \left\{ 1 + \sum_{n \in \N^*} a_n | \forall n \in \N^* , a_n \in A_n \right\} $$
is a Fr\"olicher Lie group, with regular  Fr\"olicher Lie algebra
$$\A= \left\{ \sum_{n \in \N^*} a_n | \forall n \in \N^* , a_n \in A_n \right\}.$$
Moreover, the exponential map defines a bijection $\A \rightarrow 1+\A.$  
\end{Theorem}

\noindent
\textbf{Notation:} for each $u \in \A,$ we note by $[u]_n$ the $A_n$-component of $u.$

\vskip 12pt
\noindent
\textbf{Proof.}
Let us first remark that if 
$u = \sum_{n \geq n_0} a_n \in \A$ with $n_0 \geq 1,$
then the power $u^k$ exists in $\A$ for $k \in \N^*$ since for each 
$n \in \N^*,$  the component $[u^k]_n$ is a finite sum of products of length $\leq k.$
Then, setting 
$$ u = \sum_{n \in \N^*} a_n,$$
 $1 + u \in 1 + \A$ has the inverse
$$ (1+u)^{-1} = 1 + \sum_{k \in \N^*} (-1)^k u^k \in 1 + \A.$$
The multiplication is smooth, which shows that the inversion is smooth. 
The same way, we can define an exponential
$$ exp ( u ) = 1+ \sum_{k \in \N^*} \frac{u^k}{k!},$$
a logarithm 
$$ log(1 + u) = \sum_{k \in \N^*} \frac{(-1)^{k+1}}{k} u^k$$
and fractionnal powers
$$ (1 + u)^{\frac{1}{p}} = exp\left( \frac{log(1+u)}{p} \right).$$
Now, let us turn to the regular Lie group property. 
Let $v \in C^\infty([0,1];\A).$  Let $s \in [0;1]$ and let $j = \lfloor ns \rfloor.$ We define
$$u_n(s) = \left(1+\left(s - \frac{j}{n}\right)v\left(\frac{j}{n}\right)\right) 
\prod_{i = 1}^{j}  \left(1+\frac{1}{n} v\left(\frac{j-i}{n}\right)\right).$$

We have that 
$$\lim_{n \rightarrow +\infty}\partial_s u_n(s) . u_n^{-1}(s) = \lim_{n \rightarrow +\infty} v\left(\frac{j}{n}\right)\left(1+\left(s - \frac{j}{n}\right)v\left(\frac{j}{n}\right)\right)^{-1} = v(s). $$
Moreover, the $A_m$ component of the product converges to a sum if integrals of the type  
$$\int_{1\geq s_1 \geq ... \geq s_k\geq 0} \left[ \prod_{i = 1}^k v(s_i)\right]_m (ds)^k
$$
for $k \leq m,$ which shows the convergence to a path $u$ satisfying 
$$ \partial_s u(s).u^{-1}(s) = v(s).$$ \qed

\vskip 12pt
\noindent
\textbf{Example: application to the representation of a Lie algebra.} 

We now apply this construction to a representation of a Lie algebra $\mathfrak{g}.$ For this we assume that
there exists a regular Fr\"olicher algebra $A$ (regular with respect to the abelian law $+$) such that $\mathfrak{g} \subset A$ as vector spaces,  Fr\"olicher spaces and a Lie algebras. Let us now give he construction:

\begin{itemize}

\item Let $\mathfrak{g}_1 = \mathfrak{g}$ and by induction on $n \in \N^*$ we define $\mathfrak{g}_{n+1} = [\mathfrak{g}_n,\mathfrak{g}].$ We set $\mathcal{G} = \bigoplus_{n \in \N^*} \mathfrak{g}_n.$

\item Let us fix $\forall n \in \N^*, A_{n} = A,$ and set $\A = A^{\N^*}.$
\end{itemize}

Clearly, $\mathcal{G} \subset \A$ and by Theorem \ref{regulardeformation}, the algebra $\A$ is the Lie algebra of the regular Lie group $1 + \A,$ so that, if $\mathcal{G}_1$ is the Lie subalgebra of $\A$ $Ad-$spanned by $\mathcal{G},$
$\mathcal{G}_1$ is enlargeable into a Lie subgroup of $1+\A.$

\vskip 12pt
\noindent
\textbf{Example: Formal deformation of the space of pseudo-differential operators of positive order.} An exposition  of basic facts on pseudo-differential operators can
be found in \cite{Gil}.  Let $E$ be a smooth vector bundle over a
compact manifold without boundary M. We denote by  $ PDO (M, E) $
( resp. $ PDO^k (M, E) $, resp. $ Cl(M, E) $, resp.  $ Cl^k (M, E)
$) the space of
pseudo-differential operators (resp. pseudo-differential operators
of order k, resp. classical pseudo-differential operators, resp.
classical pseudo-differential operators of order k) acting on smooth
sections of $E$. We denote by $Cl^*(M,\K^n)$,
 $Cl^{0,*}(M,\K^n)$ the groups of
the units of the algebras $Cl(M,\K^n)$ and
$Cl^{0}(M,\K^n)$.
Notice that $Cl^{0,*}(M,\K^n)$ is a CBH Lie group, and belong 
to a wider class of such groups that is studied in \cite{Glo}.

We also denote by
$Cl^{-(p),*}(M,\K^n)$ the group of invertible pseudo-differential 
operators of the
type $Id + A$, where $A \in Cl^{-p}(M,\K^n)$.
Notice that, here, the notation $Cl^{-(p),*}$ could appear 
misleading. This is 
why we feel the need to precise that this is \textbf{not} the
groups of the units of a unital algebra: this is only a regular
Lie group,
with Lie algebra $Cl^{-p}(M,\K^n)$.
A well-known fact is the following:
\begin{Lemma}
Let $v(t) \in C^\infty([0,1], Cl(M,\K^n))$ such that the order of $v(t)$ is positive for some $t\in [0;1].$ There is no smooth path $g(t) \in C^\infty( [0,1],Cl^*(M;E))$ such that $\d_tg.g^{-1}=v.$
\end{Lemma}

\vskip 12pt
\noindent
\textbf{Sketch of proof.} Assume that there is such an integral curve $g.$ Since $v$ is smooth, the order of $v$ is positive for some neighborhood of $t,$ but by the topology of $Cl^*(M,E),$ the order of $\partial_t g$ is greater than the order of $g$ only on a finite set. \qed

\vskip 12pt

In order to be complete, we have to say that particular Lie algebras of pseudo-differential operators of order $1$ are enlargeable into Lie groups of Fourier Integral operators, see e.g. \cite{Om} and the references therein, but there is no known enlargeability for $Cl^*(M,E).$ 
We propose here a formal deformation of $Cl(M,E).$
\begin{Definition}
Let $q$ be a formal parameter. 
We define the algebra of formal series 
$$Cl_q(M,E) = \left\{ \sum_{t \in \N^*} q^k a_k | \forall k \in \N^*, a_k \in Cl^k(M,E) \right\}.$$
\end{Definition}
This is obviously an algebra, graded by the order (the valuation) into the variable  $q.$ Thus, setting
$$ \A_n = \left\{ q^n a_n | a_n \in Cl^n(M,E)\right\} ,$$
we can set $\A = Cl_q(M,E)$ and state the followinc consequence of Theorem \ref{regulardeformation}:

\begin{Cor}
The group $1 + Cl_q(M,E)$ is a regular Fr\"olicher Lie group with regular 
Fr\"olicher Lie algebra $Cl_q(M,E).$
\end{Cor}

We have to say that this result is mostly inspired by the original idea of formal deformation of formal Adler series in \cite{M1,M2,M3} which will be discussed thereafter. 
Here is another result needed in the sequel. 

\begin{Theorem}\label{exactsequence}
Let $$i : K \rightarrow G$$ and $$p:G \rightarrow  H$$ such that
$$ 1 \rightarrow K \rightarrow G \rightarrow  H \rightarrow 1 $$
be an exact sequence of Fr\"olicher Lie groups, such that there is a smooth section $s : H \rightarrow G,$ and such that 
the trace diffeology from $G$ on $i(K)$ coindides with the push-forward diffeology from $K$ to $i(K).$
We consider also the corresponding sequence of Lie algebras
$$ 0 \rightarrow \mathfrak{k} \rightarrow \mathfrak{g} \rightarrow  
\mathfrak{h} \rightarrow 0 . $$
Then, 
\begin{itemize}
\item The Lie algebras $\mathfrak{k}$ and $\mathfrak{h}$ are regular if and only if the
Lie algebra $\mathfrak{g}$ is regular;
\item The Fr\"olicher Lie groups $K$ and $H$ are regular if and only if the Fr\"olicher Lie group $G$ is regular.
\end{itemize}

\end{Theorem}

\vskip 12pt
\noindent
\textbf{Proof.}

$ \bullet$ \underline{At the level of Lie algebras}, 
we set $i' = D_ei$ and $p' = D_ep$. We have 
$$i'(\mathfrak{k}) = Ker\left(s' \circ p'\right)$$ 
so that $$ \mathfrak{g} \sim \mathfrak{k} \oplus \mathfrak{h}$$
which proves the first equivalence (through integration componentwise).

$\bullet$ \underline{At the level of Lie groups} We follow the proof given in \cite{KM}, section 38.6 in the case of Lie groups in the convenient setting. 

- If $G$ is regular, let $Y \in C^\infty(\R ; \mathfrak{h})$ then there exists $y_0 \in C^\infty(\R , \mathfrak{g})$ such that $$\partial_t y_0(t) . y_0(t)^{-1} = s'(Y(t)).$$
Then, $$ y(t) = p \circ y_0 (t)\in C^\infty(\R, H)$$ integrates $Y.$ Moreover, 
let $Z \in C^\infty(\R,\mathfrak{k}).$ then $i'(Z)$ is integrated into a smooth curve $i(z) \in C^\infty(\R, i(K)).$ Since the trace diffeology from $G$ on $i(Z)$ coincides with the push-forward diffeology from $K$ to $i(K),$ the path $z$ is smooth on $K$ and integrates $Z.$ 

- If $H$ and $K$ are regular, let $X \in C^\infty(\R, \mathfrak{g}).$
Let 
$$ Y = p'(X) \in C^\infty(\R,\mathfrak{h})$$
and let $h\in C^\infty(\R, H)$ integrating $Y$ on $H.$ 
Let $Z \in C^\infty(\R, \mathfrak{k})$ such that 
$$ i'(Z)(t) = Ad_{s \circ h(t)^{-1}} \left( X(t) - D_{h(t)} s \left(Y(t).h(t)\right).(s\circ h)(t)^{-1}\right) .$$
We set $k \in C^\infty(\R, K)$ that integrates $Z$ on $K$ 
and define $$g(t) = (s \circ h)(t) . (i \circ k)(t)\in C^\infty(\R,G).$$
Then, $$\partial_tg(t). g(t)^{-1} = X(t)$$
which shows that $G$ is regular.

\section{Path liftings} \label{3}

\subsection{The general setting for path-lifting} \label{3.1}

Let $X$ and $M$ be two differentiable spaces,
and $\pi: X \rightarrow M$ a smooth surjective map. 
If $\gamma$ is a path, 
we set $\gamma^{-1}(t) = \gamma(1-t).$
We define 
$\B(X;M)$  the set of couples of the type  $(\gamma,x)$ such that $\gamma$ is a smooth 
path on $M$ and $x \in \pi^{-1}(\gamma(0)).$ It is a natural (trace) differentiable space. 
  
\begin{Definition} \label{pl},
A \textbf{path-lifting} $L: (\gamma;x) \mapsto L_x(\gamma)$ is a smooth map
from $\B(X;M)$ to the set of smooth paths on $X,$ which satisfy the following properties :
\vskip 8pt

(i) $\forall \gamma \in C^\infty([0,1],M),$ $\forall x \in \pi^{-1}(\gamma(0)),$ 
$$ \pi \circ (L_x(\gamma)) = \gamma;$$

(ii) $\forall \gamma \in C^\infty([0,1],M)$ such that $\gamma = \gamma' \vee \gamma'',$
$\forall x \in \pi^{-1}(\gamma(0)),$ 
 $$ L_x(\gamma) = L_{L_x(\gamma'')(1)}(\gamma') \vee L_x(\gamma''),$$ 

(iii) $\forall \gamma \in C^\infty([O,1],M),$ $\forall g \in C^\infty([0,1], [0,1])$
such that $g$ is monotone, $\forall x \in \pi^{-1}(\gamma(0)),$ 
$$ L_x(\gamma \circ g) = L_x(\gamma) \circ g.$$ 

(iv) $\forall (\gamma, x) \in \B(X,M),$
$$L_x(\gamma^{-1} \vee \gamma) (1) = L_x(\gamma)(0)=x.$$

(v) $\forall (\gamma,x)\in \B(X,M),$ $L_x(\gamma)(0) =x.$

(vi) Let $\gamma \in C^\infty([0;1];M)$. 
$$(\exists x \in \pi^{-1}(\gamma(0)) \hbox{, } L_x(\gamma)(1) = x) \Leftrightarrow
  (\forall x \in \pi^{-1}(\gamma(0)) \hbox{, } L_x(\gamma)(1) = x)$$
\end{Definition}

\vskip 8pt
One can recognize here all the basic key properties of the horizontal lift
of a connection on a fiber bundle when the fiber is a compact manifold without boundary \cite{KMS}.
We need to precise an ambiguity due to the ``smooth'' setting: the product of
paths is often developed for \textbf{continuous} (and not smooth) paths, because the condition
$\gamma(1) = \gamma'(0)$ is not sufficient to get a smooth path $\gamma' \vee \gamma$ if $\gamma$ and 
$\gamma'$ are smooth. This is why we need to reparametrize paths into paths 
which are stationary at each endpoint. This does not change anything 
for path lifting (except the parametrization), because
of property (iii). For the details in a context of 
connections on principal bundles, 
see e.g. \cite{Ma}.  

We now have to extend the classical 
construction of the holonomy group of a connection
to the context of path liftings.

Let $x \in X,$ $m = \pi(x).$ We set 
\begin{eqnarray*} \lo(m,M) & = & \{\gamma : [0;1] \rightarrow M \hbox{ such that } \gamma
\hbox{is a loop based at {\it m}, stationary at endpoints} \} \end{eqnarray*}
\vskip 8pt
 
Let $\lo_0^L(m;M) = \{ \gamma \in \lo(m;M) \hbox{ such that } L_x(\gamma)(1) = x \}.$
Notice that all the loops of the type $\gamma^{-1} \vee \gamma$ are in $\lo_0^L(m;M)$ 
by (iv) of Definition \ref{pl}.
\begin{Definition} For $\gamma, \gamma'\in \lo(m;M),$
\begin{eqnarray*} \gamma \sim \gamma' & \Leftrightarrow & \exists x \in \pi^{-1}(m), L_x(\gamma)(1) = L_x(\gamma')(1) 
\end{eqnarray*} 
\end{Definition}

Using (ii) of Definition \ref{pl}, we get immediately: 
\begin{Proposition}
\begin{eqnarray*} \gamma \sim \gamma' & \Leftrightarrow & (\gamma')^{-1} \vee \gamma \in \lo_0^L(m;M) \hbox{ . }\end{eqnarray*}
\end{Proposition}

Using (vi) of Definition \ref{pl}, we see: 
\begin{Proposition}
\begin{eqnarray*} \gamma \sim \gamma' & \Leftrightarrow & \forall x \in \pi^{-1}(m), L_x(\gamma)(1) = L_x(\gamma')(1) \hbox{ . }\end{eqnarray*}
\end{Proposition} 

Finally, by (iii) of Definition \ref{pl}, we get the last easy property:

\begin{Proposition} if $\tilde \gamma$ is a reparametrization of $\gamma \in \lo(m;M),$
$\tilde \gamma \sim \gamma.$ \end{Proposition}

Noticing that this relation is symmetric and transitive,  
we set
$$\Hol_x^L = \lo(m;M) / \sim.$$

\vskip 8pt
$\bullet$ \textit{push-forward of $\vee$ and of the inversion of paths}

First, we notice that, 

- if $\gamma \sim \gamma'$, $L_{L_x(\gamma)(1)}(\gamma^{-1})(1)= x = L_{L_x(\gamma')(1)}(\gamma'^{-1})(1)$ 
hence $\gamma^{-1} \sim \gamma'^{-1}.$

- if $\gamma \sim \gamma'$ and $\delta \sim \delta',$  $\forall x \in \pi^{-1}(m),$
$$
L_x(\gamma \vee \delta)(1) =  L_{L_x(\delta)(1)}(\gamma)(1) =  
L_{L_x(\delta')(1)}(\gamma)(1) =  L_{L_x(\delta')(1)}(\gamma')(1)
   =  L_x(\gamma' \vee \delta')(1). $$
Hence $\gamma \vee \delta \sim \gamma' \vee \delta'.$

\vskip 8pt 
$\bullet$ \textit{Associativity.} Given three loops $\gamma,$ $\gamma'$ and $\gamma''$, 
$\gamma \vee (\gamma' \vee \gamma'')$ and $(\gamma \vee \gamma') \vee \gamma''$ differ by parametrizations.
Hence,  $\gamma \vee (\gamma' \vee \gamma'') \sim (\gamma \vee \gamma') \vee \gamma''.$ 

\vskip 8pt
We note also by $\vee$ and ${}^{-1}$ the push forward of  $\vee$ and ${}^{-1}$ onto $\Hol^L_x.$ We can now state, by push-forward of the differentiable structure of $L(M,m):$ 
\begin{Theorem}  $(\Hol_x^L, \vee)$ is a diffeological group, and the inversion ${}^{-1}$ is differentiable.
\end{Theorem}    
 
\subsection{Comments and remarks}

\vskip 5pt
\noindent
$\bullet$ Let $\gamma \in \lo(m;M).$ The map $L_{(.)}(\gamma)$ is a smooth map on 
$\pi^{-1}(m)$, with smooth inverse $L_{(.)}(\gamma^{-1}).$ Then, for fixed 
$m \in M$, $L$ defines a map $ \lo(m;M) \rightarrow Diff(\pi^{-1}(m)),$ where 
$Diff(\pi^{-1}(m))$ is the group of (diffeological) diffeomorphisms of the fiber.
Moreover, 
\begin{Proposition} Let $m \in M$ and $\gamma \in \lo(m;M).$
$$L_{(.)}(\gamma' \vee \gamma) = L_{(.)}(\gamma') \circ L_{(.)}(\gamma),$$
\end{Proposition}

\noindent 
\textbf{Proof.} Straightforward from $\forall \gamma \in C^\infty([0,1],M)$ such that $\gamma = \gamma' \vee \gamma'',$
$\forall x \in \pi^{-1}(\gamma(0)),$ 
 $$ L_x(\gamma) = L_{L_x(\gamma'')(1)}(\gamma') \vee L_x(\gamma'')\qed$$

\begin{Proposition}
Let $m \in M$ and $\gamma \in \lo(m;M).$
\begin{eqnarray*}L_{(.)}(\gamma) = Id_{\pi^{-1}(m)} & \Leftrightarrow & \gamma \in \lo_0^L(m;M) \end{eqnarray*}
\end{Proposition}

\noindent 
\textbf{Proof.} Straightforward from: Let $\gamma \in C^\infty([0;1];M)$. 
$$(\exists x \in \pi^{-1}(\gamma(0)) \hbox{ such that } L_x(\gamma)(1) = x) \Leftrightarrow
  (\forall x \in \pi^{-1}(\gamma(0)) \hbox{ such that } L_x(\gamma)(1) = x) \qed $$

Then, we have:
\begin{Theorem} Let $m \in M$ and $x \in \pi^{-1}(M).$ The map $L$ induces a quotient map $$ \Hol_x^L \rightarrow Diff(\pi^{-1}(m)).$$
\end{Theorem}

\vskip 5pt
\noindent
$\bullet$ Let $(P,M,G)$ be a (classical) finite dimensional principal bundle 
of basis $M$ and of structure group $G.$ Any connection on $P$ induces a 
path-lifting. 
It seems that the inverse induction is not elementary, and is not true in general,
because the structure group $G$ (which models the fiber) 
is viewed here only as a manifold, 
without any group structure. Maybe a stronger analysis could give more details on the correspondence, 
for example up to homotopy, between general path liftings and path liftings induced by connections.

\vskip 5pt
\noindent
$\bullet$ Let $(P,M,N)$ be a finite dimensional fiber bundle of basis $M$ 
and of typical fiber $N$. Then, if $Diff(N)$ is a Lie group, 
$P \times_{N} Diff(N)$ is a principal bundle and there is a bijection between
fiber bundle connections on $(P,M,N)$ and (classical) connections on 
$P \times_N Diff(N).$
This bijection is established in e.g. \cite{KM} or \cite{KMS}.
But, if $N$ is not compact, 
horizontal lifts of paths are not well-defined for an arbitrary connection,
and one has problems to define a holonomy group where as curvature 
elements are well-defined. 
This comes from the definition of horizontal lifts, 
which are defined as solutions of a differential equation: 
$ \tilde \gamma \in C^\infty([0;1];P)$ is a horizontal lift of 
$\gamma \in C^\infty([0;1];M)$ if $\pi(\tilde \gamma) = \gamma$ 
and if $D\tilde \gamma$ is horizontal. With such a definition, it is obvious 
that two horizontal lifts of a same path $\gamma$ differ by their starting point, 
but also that all starting points are not good to define horizontal 
lifts for $\gamma$ when $N$ is not compact. 
With our setting, such problems are avoided since only ``good'' connections, 
that is connections for which horizontal lifts exist, are considered.
     
\vskip 5pt
\noindent
$\bullet$ When $X$, a diffeological space, is equipped with a relation of equivalence $\rel$, the
quotient space $M = X / \rel$ (with quotient projection $\pi : X \rightarrow M$) 
is also a diffeological space by push-forward of the diffeology of $X.$
Then, our notion of path-lifting is also valid here, as well as holonomy. In a future work, we shall precise the role of the condition (vi) on the existence of an isomorphism
between two equivalence classes $\pi^{-1}(m)$ and $\pi^{-1}(m')$ for $(m,m') \in M^2.$ A counterexample to invariance of the holonomy group will be developped elsewhere, in the framework of infinite configuration spaces.  

\subsection{Homotopy, fundamental group and holonomy of a path-lifting} \label{pafin}

As mentioned in \cite{Sou}, the notion of 
homotopy of paths can be adapted straightway 
from the category of topological spaces to the category of diffeological spaces.
These two notions coincide on the subcategory of 
smooth finite dimensional (paracompact)
manifolds, since smooth or continuous homotopy gives 
the same equivalence classes of maps
in this restricted class of objects.
We also recall that one can adapt straightway the
definition of arcwise connected components to the setting of
diffeological spaces:

\begin{Definition}
Let $(X, \p)$ be a diffeological group. 
Let $(x, y) \in X^2.$ $x$ and $y$ are in the same (arcwise) connected component
if there is a smooth path $\gamma$ on $X$ starting from $x$ and ending on $y.$
\end{Definition}
 
Let us now recall the definitions and the key properties of the homotopy of loops
in the context of diffeological spaces:

\vskip 12pt
\noindent
\begin{it}
Let $X$ be a diffeological space.
Let $\gamma$ and $\gamma'$ be two smooth loops on $X$, based on $x$. 

$\bullet$ A \textbf{homotopy} between $\gamma$ and $\gamma'$ is a smooth map 
$$ H : S^1 \times [0;1] \rightarrow X$$
such that $H(1,.) = x,$  $H(.,0) = \gamma$ and $H(.,1) = \gamma'.$

$\bullet$ $\gamma$ and $\gamma'$ are called homotopic if there exists a homotopy between $\gamma$ and $\gamma'.$

$\bullet$ $\pi_1(X,x)$ is the set of connected components of smooth loops.
This is the space of equivalence classes of loops modulo homotopy. If we consider only
loops constant at endpoints, $\pi_1(X,x)$ gets
a group structure
induced by the composition of paths $\vee.$

$\bullet$ $\pi_1(X,x)$ can be identified to the set of connected components of 
${\lo}(X,x).$

\end{it}

\vskip 12pt

From this last property, with the notations of section \ref{3.1}, 
any path lifting $L$ induces a (onto) map from $\pi_1(M,m)$
to the connected components of $ \mathcal{H}^L_x.$ 

\begin{Definition}
A path lifting is \textbf{flat} at $x \in X$ if the connected components 
of $\mathcal{H}^L_x$ are made of singletons.

A path lifting in \textbf{totally flat} at $x \in X$ if the set $\mathcal{H}^L_x$ 
has only one element.
\end{Definition}

With this definition, we get the following obvious statement:

\begin{Proposition}
If $L$ is flat at $x$, there is a surjective group morphism
$$  \pi_1(M,m) \rightarrow \mathcal{H}^L_x$$
induced by the path lifting $$L : \mathcal{L}(M,m) \rightarrow \hbox{ paths starting at }x.$$ 
\end{Proposition}

\section{Diffeological principal bundles with regular (Frölicher) groups}
\label{conn}
Let $P$ be a diffeological space and  let $G$
be a regular Frölicher Lie group, 
with a differentiable right-action 
$P\times G \rightarrow P,$ such that 
$\forall (p,p',g)\in P\times P \times G,$ we have
$p.g=p'.g \Rightarrow p=p'$.
Let $M=P/G,$ equipped with the quotient diffeology.

\begin{Proposition}
Let $V$ be a vector space. $G$ acts smoothly on the right on $\Omega(P,V)$ setting
$$ \forall (g,\alpha) \in  \Omega^n(P,V) \times G, \forall p\in \p(P), (g_*\alpha)_{g.p} = \alpha_p \circ (dg^{-1})^n .$$
\end{Proposition}
\noindent
\textbf{Proof.}
$G$ acts smoothly on $P$ so that, if $p \in \p(P),$ $g.p \in \p(P)$.
The right action is now well-defined, and the smoothness is trivial. 

\begin{Definition}
Let $\alpha \in \Omega(P;\mathfrak{g}).$ The differential form $\alpha$ is \textbf{right-invariant} if and only if, for each $p \in \p(P),$ and for each $g \in G,$
 $$\alpha_{g.p} = Ad_{g^{-1}} \circ g_*\alpha_p.$$
\end{Definition}

 Let $p \in P$ and 
$\gamma$ a smooth path in $P$ starting at $p.$

\begin{Definition}
A \textbf{connection} on $P$ is a $\mathfrak{g}-$valued $1-$form $\theta,$ right-invariant, such that, for each $ g \in \mathfrak{g},$ for any path $c : \R \rightarrow G$ such that $\left\{ \begin{array}{ccc} c(0)& = & e_G \\ c'(0)&= &g \end{array} \right. ,$ and for each $p \in P,$ $$\theta((p.c(t))'_{t = 0})=g.$$ 
\end{Definition}
Now, let us turn to holonomy. Let $p \in P$ and 
$\gamma$ a smooth path in $P$ starting at $p,$ defined on $[0;1].$
Let $H\gamma (t) = \gamma(t)g(t)$ where $g(t) \in C^\infty([0;1];\mathfrak{g})$ is a path satisfying the differential equation: 
$$\left\{ \begin{array}{c} \theta \left( (H\gamma(t))' \right) = 0  \\ H\gamma(0)=\gamma(0) \end{array} \right.$$
The first line of the definition is equivalent to the differential equation 
$g^{-1}(t)(g(t))' = -\theta(\gamma ' (t))$, and the second to $g(0)=e_G,$ which is integrable. 
This shows that horizontal lifts are well-defined, as in the case of manifolds. 
Moreover, the map $H(.)$ defines trivially a path-lifting. This enables us to consider the 
holonomy group of the connection.
Notice that a straightforward adaptation of the arguments of \cite{Ma} shows that the 
holonomy group is invariant (up to conjugation) under the choice of the basepoint $p.$
We note it now $\Hol,$ omitting the basepoint $p$ and the connection $\theta$ in our notations since 
they are assumed fixed.
 
\section{Curvature and the Lie algebra of the Holonomy group}

Now, we assume that $dim(P)\geq 2.$ We fix a connection $\theta$ on $P.$

\begin{Definition}
Let $\alpha \in \Omega(P;\mathfrak{g})$ be a $G-$invariant form. Let $\nabla \alpha = d\alpha - {\frac{1}{ 2}}[\theta,\alpha]$ be the horizontal derivative of $\alpha.$
We set $$ \Omega = \nabla \theta $$ the curvature of $\theta.$ 
\end{Definition}
  
\subsection{Reduction of the structure group}
  We now turn to reduction of the structure group, adapting a theorem from \cite{Ma2}:
 
 \begin{Theorem} \label{Courbure}
We assume that $G_1$ and $G$ are regular Fr\"olicher groups 
with regular Lie algebras $\mathfrak{g}_1$ and $\mathfrak{g}.$
Let $\rho: G_1 \mapsto G$ be an injective morphism of Lie groups.
If there exists a connection $\theta$ on $P$, with curvature $\Omega$, such that, for any smooth 1-parameter family $Hc_t$ of horizontal paths starting at $p$, for any smooth vector fields $X,Y$ in $M$,  
\begin{eqnarray} s, t \in [0,1]^2 & \rightarrow & \Omega_{Hc_t(s)}(X,Y)  \label{g1}\end{eqnarray} 
is a smooth $\mathfrak g_1$-valued map (for the $\mathfrak g _1 -$ diffeology),
\noindent
and if $M$ is simply connected, then the structure group $G$ of $P$ reduces to $G_1,$ and the connection $\theta$ also reduces. 
\end{Theorem}  
  Before giving the proof of this theorem, we need three lemmas, 
which are well-known results in finite dimensions (see e.g. \cite{Li}). Following \cite{KM}, 
if $ C^{\infty}_x([0,1],M) $ is the set of smooth paths on $M$ starting at $x$,
we note $$Pt : C^{\infty}_x([0,1],M) \times [0,1] \times \pi^{-1}(x)) \rightarrow P$$ the parallel transport with respect to $\theta$, which is a smooth map.
Let us now describe the skeleton of the proof: the key tools for the definition of the local trivializations 
needed to reduce the principal bundle $P$ is given in lemma \ref{plaques}, 
which is inspired by the computations in \cite{KM} in the case of a vanishing curvature 
and \cite{Ma2} in the general case. Lemmas \ref{horizontal} and \ref{final} deal with local 
description of horizontal lifts. Finally, in the proof of Theorem \ref{Courbure}, we define 
a family of local trivializations of $P$ using Lemma \ref{plaques}, and check that it has the 
desired properties by Lemmas \ref{horizontal} and \ref{final}. Till the end or the proof of 
Theorem \ref{Courbure}, we use the notations defined in the beginning of this paragraph.

\begin{Lemma} \label{plaques}
Let $p \in P$. Let $x$ be the basepoint of $p$. 
Let $\varphi:U\rightarrow M \in \p$ be a plot of the diffeology of $M$ with star-shaped domain $U$, 
that we identify (for the sake of simplicity) with a  
star-shaped neighborhood of 0 in $\R ^n.$ 
Let $u \in U$ and $t \in [0,1]$. We define $f(u,t)= tu \in U$. 

Let \begin{eqnarray*}
\psi : U & \rightarrow & P \\
     u & \mapsto & Pt(f(u,.), 1  , p) .\end{eqnarray*}

Let \begin{eqnarray*} \Psi : U \times G & \rightarrow & \pi^{-1}(U)\\
(u,g) & \mapsto & \psi(u).g \quad ,\end{eqnarray*}

and $$ \tilde\Psi = \Psi \circ (Id_{U} \times \rho).$$

Then, $\psi$ is a plot of the diffeology of $P.$ Moreover, $\theta \circ D\psi$ is a smooth 
$\mathfrak g_1$-valued form on $U$. 
\end{Lemma}
We need now to know how horizontal lifts of paths behave in the diffeology of P, and more precisely in with respect to the structure of $G_1.$

\begin{Lemma}\label{horizontal}
We assume that $U$ is convex.
\begin{item}
(i) Given a path $\alpha: [0,1] \rightarrow U$ starting at $x$, if $H\alpha$ is its horizintal lift starting at $p$, we have $H\alpha(1) \in \Psi ( U \times \rho(G_1))$, and there exists a smooth path $H\alpha_1 :[0,1] \rightarrow U \times G_1$ such that $H\alpha=\tilde\Psi  \circ H\alpha_1.$
\end{item}
\begin{item}
(ii) Let $h:[0,1]^2 \rightarrow U $ be an homotopy equivalence between two paths $h(0,.)$ and $h(1,.)$ starting at $x$ and finishing in $U$. Let $Hh(0,.)$ and $Hh(1,.)$ be their horizontal lifts starting at $p$. Then, there is $g_1 \in G_1$ such that $\tilde\Psi^{-1}(Hh(0,1)) = \tilde\Psi^{-1}(Hh(1,1)).\rho(g_1).$    
\end{item}
\end{Lemma}

Then, the following lemma will be useful when dealing with homotopy:

\begin{Lemma}\label{final}
Let $\alpha$ and $\beta$ be two paths on $U$. Let $q_\alpha \in \pi^{-1}(\alpha(0))$ and $q_\beta \in \pi^{-1}(\beta(0))$. Let $H\alpha$ and $H\beta$ be the horizotal lifts of $\alpha$ and $\beta$ starting at  $q_\alpha$ and $q_\beta$. We set $\tilde \Psi^{-1}\circ H\alpha = (\alpha,\gamma_\alpha)$ and $\Psi^{-1}\circ H\beta = (\beta,\gamma_\beta)$. 

Let $g = \gamma_\beta^{-1}(0).\gamma_\alpha(0)$, with $^{-1}$ as 
the inverse map in $G$. Then, for any $t \in [0,1]$, there exists $g_1(t), g'_1 \in G_1$ such that 
$\gamma_\beta(t) = \gamma_\alpha(t).g_1^{-1}(t).g.g_1'(t)$. Moreover, 
the maps $t \mapsto g_1(t)$ and $t \mapsto g'_1(t)$ are smooth in $G_1$.
\end{Lemma}
\noindent 
Let us now give the proofs of the three lemmas, and then the proof of Theorem \ref{Courbure}:
\vskip 10pt
\noindent
\textbf{Proof of Lemma \ref{plaques}:}

We already know that $\Psi: U \times G \rightarrow P$ is a smooth map, since $Pt$ is smooth.   

Let us calculate $D\Psi^{-1} \circ \theta \circ D\psi$.
Let $c:]-\epsilon, \epsilon[ \rightarrow U$ be a smooth path such that $c(0) = u \in U$. Let $h(t,s) = f(c(s),t)$. 
Let $\tilde \theta$ be the pull-back of $\theta$ by $\Psi$ on $U$.

Then, following the proof of the claim of \cite{KM}, theorem 39.2, with a connection with non vanishing curvature, see e.g. \cite{Ma2}, 
we have:

\begin{eqnarray*} \partial_s(h^*\tilde \theta)(\partial_t) & = &  \partial_t(h^*\tilde \theta)(\partial_s) - d(h^*\tilde \theta)(\partial_t,\partial_s) - (h^*\tilde \theta)([\partial_t,\partial_s]) \\
& = & \partial_t(h^*\tilde \theta)(\partial_s) - d(h^*\tilde \theta)(\partial_t,\partial_s)\\
& = & \partial_t(h^*\tilde \theta)(\partial_s) + ad_{(h^*\tilde \theta)(\partial_t)}((h^*\tilde \theta)(\partial_s)) - (\Psi^*\Omega)(h_*(\partial_s),h_*(\partial_t)) \\
& = & \partial_t(h^*\tilde \theta)(\partial_s) - (\Psi^*\Omega)(h_*(\partial_s),h_*(\partial_t)) \quad \hbox{since } (h^*\tilde \theta)(\partial_t)= 0 .\end{eqnarray*}

Thus, 

$$\partial_s(\Psi \circ \psi \circ c) = (\partial_s c(s), \partial_s \tilde\gamma(1,s)),$$
remarking that $(u,e) = \Psi^{-1} \circ \psi(u)$. 
We now calculate $\partial_s \tilde\gamma(1,s)$, 
\begin{eqnarray*}
\partial_s \tilde\gamma(1,s) & = & \int_0^1 \left( \partial_s(h^*\tilde \theta)(\partial_t) \right)(t
) dt \\
& = & \int_0^1 \Big(\partial_t(h^*\tilde \theta)(\partial_s) 
- (\Psi^*\Omega)(h_*(\partial_s),h_*(\partial_t))\Big)(t) dt
\\
& = & (h^*\tilde \theta)(\partial_s)(1,s) - \int_0^1 (\Psi^*\Omega)(h_*(\partial_s),h_*(\partial_t))(t) dt.\end{eqnarray*}
Finally, we have: 
\begin{eqnarray} \theta(\partial_s(\psi \circ c)) &=& (h^* \tilde \theta)(\partial_s h(1,s),\partial_s \gamma(1,s))\\
& = & \int_0^1 \Big(\Omega(h_*(\partial_s),h_*(\partial_t))(t) \Big)dt. \label{evolution}
 \end{eqnarray}
Since $\mathfrak g_1$ is complete, this integral exists and belongs to $\mathfrak g_1$.
\qed
\vskip 10pt
\noindent
\textbf{Proof of the Lemma \ref{horizontal}:} 

\textit{(i)} We have that $\theta(\partial_s(\psi \circ \alpha))$ is an integral on the curvature elements (see the proof of the last lemma). 
Looking at this result more precisely, reparametrizing equation (\ref{evolution}), setting $c=\alpha$, we have that
\begin{eqnarray}\theta(\partial_s(h(s,t))) & = & \int_0^t \Big(\Omega(h_*(\partial_s),h_*(\partial_t))(u) \Big)du,\end{eqnarray} and hence that $$\partial_t\Big(\theta(\partial_s(h(s,t))\Big) = \Omega(h_*(\partial_s),h_*(\partial_t)).$$
Recall that $ \rho^{-1}\circ\Omega(h_*(\partial_s),h_*(\partial_t)) $ is smooth.
Integrating this equality in $G_1$ instead of $G$, we get a path $\alpha_1$ in $U \times G_1$. Then we consider the following differential equation, that defines $H\alpha_1$: 
$$ \left\{ \begin{array}{ccl}H\alpha_1(0)&=&e\\
 H\alpha_1(t) & = & (\alpha(t),\gamma(t)) \in U \times G_1\\
D\rho\big( \gamma(t)^{-1}\partial_t\gamma(t) \big) &=& Ad_{(\rho \circ \gamma^{-1})(t)} \Big( \theta \circ D\Psi\circ(Id \times\rho)(\partial_t\alpha_1)(t)\Big) , \end{array} \right.$$
setting  $H\alpha=\Psi \circ (Id_U \times \rho) \circ H\alpha_1,$ we get (i).

\textit{(ii)} comes easily from the continuity of the horizontal lift of paths, using the fact that $U$ is contractible, and applying \textit{(i)} to the path 
\begin{equation} \label{hg1} c(t)= \left\{\begin{array}{ccl} Hh(0,3t) & \hbox{if} & t \in [0,1/3] \\
   Hh(3t-1,1)& \hbox{if} & t \in [1/3,2/3]\\
 Hh(1,3-3t)& \hbox{if} & t \in [2/3,1].\end{array}\right.\end{equation}
\qed
\vskip 10pt
\noindent
\textbf{Proof of Lemma \ref{final}:} Reparametrizing Lemma \ref{horizontal}, (i), we have that, for $t \in [0,1]$, there exists $g_1(t), g'_1(t) \in G_1$ such that $\gamma_\alpha(t) = \gamma_\alpha(0).g_1(t)$ and $\gamma_\beta(t) = \gamma_\beta(0).g'_1(t)$. Then, we get $\gamma_\beta(t) = \gamma_\alpha(t).g_1^{-1}(t).g.g_1'(t)$. By equation (\ref{hg1}), we have that the paths $t \mapsto g_1(t)$ and $t \mapsto g'_1(t)$ are smooth paths with values in $G_1$. \qed   
\vskip 10pt
\noindent
\textbf{Proof of Theorem \ref{Courbure}:}
Let $p_o \in P,$ and 
let $X$ be a family of paths $\gamma_x$  in $M$, indexed by $ x \in M,$ 
starting at $x_0 = \pi(p_0)$ and ending at $x.$ 
Let $X_1 = \bigcup_{x \in M} H\gamma_x(1)$ and 
let $P_1 = X_1.\rho(G_1).$ Let $\delta$ be an arbitrary path starting at 
$x_0$, and let $x_1 = \delta(1).$ The path $\gamma_{x_1}\vee \delta^{-1}$ is 
null-homotopic since $M$ is simply connected. So that, using Lemma \ref{horizontal},
 there exists $g_1 \in G_1$ such that 
$\gamma_{x_1}\vee \delta^{-1}(1)=p_0.\rho(g_1).$ So that $\delta(1)= \gamma_{x_1}(1).g_1^{-1}$ 
and $\pi^{-1}(x_1) = p_0.G_1.$
Moreover, by Lemma \ref{plaques}, the maps $\varphi: U \rightarrow P$ are $P_1-$valued. We define
a smooth diffeology on $P_1$ which is generated by the push-forward diffeologies of 
the subsets of the type $Im(\varphi).G_1,$ induced by the maps $(u,g_1)\in U \times G_1 \mapsto \varphi(u).G_1.$ 
With this diffeology:

$\bullet$ the inclusion map $P_1 \rightarrow P$ is smooth

$\bullet$ the horizontal lift map $\alpha \mapsto H\alpha$ is a smooth map by trivial application of 
Lemma \ref{final}

$\bullet$ The connection $\theta$ restricts to a smooth $\mathfrak{g}_1$-valued form by Lemma \ref{plaques}.
This ends the proof of the reduction theorem, since $G_1$-right invariant trivially. \qed

\subsection{Ambrose-Singer theorem in the Fr\"olicher setting}

We can now state the announced Ambrose-Singer theorem. In this theorem, we review the key results of this paper:

\begin{Theorem}
\label{Ambrose-Singer}
Let $P$ be a principal bundle of basis $M$ with regular Fr\"olicher structure group $G$ with regular Lie algebra $\mathfrak{g}.$
Let $\theta$ be a connection on $P$ and $L$ the associated path-lifting.

    \begin{enumerate}
\item For each $p \in P,$ the holonomy group $\Hol_p^L$ is a 
diffeological subgroup of $G$, which does not depend on the choice of 
$p$ up to conjugation.

 \item There exists a second holonomy group $H^{red},$ $\Hol \subset H^{red},$ 
which is the smaller structure group for which there is a subbundle $P'$ to 
which $\theta$ reduces. Its Lie algebra is spanned by the curvature elements, i.e.
it is the smallest integrable Lie algebra which contains the curvature elements. 

\item If $G$ is a Lie group (in the classical sense) of type I or II in 
the terminology of Robart \cite{Rob}, 
there is closed Lie subgroupgroup $\bar{H}^{red}$ (in the classical sense) such that $H^{red}\subset \bar{H}^{red},$ 
whose Lie algebra is the closure in $\mathfrak{g}$ of the Lie algebra of $H^{red},$ 
which is the smaller closed Lie subgroup of $G$ among the structure groups 
of closed subbundles $\bar{P}'$ of $P$ to which $\theta$ reduces.
\end{enumerate}
\end{Theorem} 

\noindent
\textbf{Proof.}

\begin{enumerate} 

\item is proved in section \ref{conn}

\item Let $P\Hol$ be the set of elements of $P$ that are joint to $p$ by a 
horizontal path. $P\Hol$ is obviously a principal bundle with structure group 
$\Hol$ (or a ``structure quantique'' using the terminology of Souriau \cite{Sou}). 
Notice that we do not assume here local trivializations on $P\Hol.$ 
By Theorem \ref{Courbure}, for each regular Fr\"olicher Lie subgroup $G_1$ of $G$, 
with Lie algebra $\mathfrak{g}_1,$ if $\mathfrak{g}_1$ is regular, 
the connection $\theta$ reduces to the bundle 
$$P\Hol \times_\Hol G_1 = (P\Hol \times G_1)/\Hol .$$
 The family $\mathfrak{G}$ of such Lie groups $G_1$ is not empty 
since $G \in \mathfrak{G},$ and it is obviously filtering for $\supset.$ 
So that $\mathfrak{G}$ has a minimal element $H^{red}$ for $\subset.$ By 
Theorem \ref{Lie}, the Lie algebra of $H^{red}$ is the 
smaller regular Lie algebra which contains the curvature elements. 

\item is a straightforward application of the arguments of \cite{Ma2}.
\end{enumerate}

\section{Application to the KP hierarchy}

In the following, $R$ is an algebra of smooth functions over a field 
$k = \R$ or $\C$ which is endowed with a Fr\"olicher structure derived 
from a Fr\'echet topology, and $\d : R \rightarrow R$ is a smooth $k-$linear 
map that satisfies the Liebnitz rule $\d ( fg) = \d (f) g + f \d (g)$. We 
require : 

(i) $R$ has a unit element,

(ii) for any $f \in R$, there is $g \in R$ such that $\d g = f$,

(iii) $R$ is closed under exponentiation.

\subsection{The Lax equation}

Let $D = R[\d ]$ be the algebra of differential operators, and let 
$E = R((\d ^{-1}))$ be the algebra of Adler series, endowed with 
the product topologies. In \cite{L}, Lax shows that a 1-parameter 
family $\{P_t\}$ of $E$ is an isospectral family if and only if $\{P_t\}$ 
satisfies the Lax equation $$ {d P_t \over dt} = [Q(t),P(t)],$$ 
where $Q(t)$ is a family of differential operators. The problem 
reduces to the following (see e.g. \cite{M2}): 

Let $E^{(-1)}$ be the subalgebra of $E$ made of series of order $-1$, i.e. 
$$ E^{(-1)} = \{ \sum_{n = 1}^{+\infty} a_{-n} \d ^{-n} | a_{-n} \in R \}$$
Then, $$ E : D \oplus E^{(-1)},$$
and any $P \in E$ decomposes as a sum $P = P_+ + P_-$, where $P_+ \in D$ and $P_- \in E^{(-1)}$. 
Let $(t_1, t_2,...)$ be a countable family of variables. 
We define $$T = \bigoplus_{n \in \N^*} k.t_i $$ the algebraic $k-$ vector space spanned by the variables
$t_1,t_2,...$ .  The Lax equation 
reduces to find 
  $$L = \d + a_{-1}\d^{-1}+ a_{-2}\d^{-2}+...$$
where each $ a_{-n} \in \mathcal{R} = R[[t_1, t_2, .... ]] = R \otimes_k T$; 
 $\mathcal{R}$ is a the algebra of formal power series with valuation defined by 
$val(t_n) = n$, and such that $$ {dL \over dt_n} = [L^n_+, L], $$
which can be also written as 
\begin{equation} \label{lax1} d L = [Z_+, L] = [Z_-,L], \end{equation}
where $Z = \sum_{n \in \N} L^n dt_n,$ $Z_+= (Z)_+$ and $Z_- = (-Z)_-.$ Moreover, if we want to get a family of isospectral deformations (which was the motivation to this setting), we have also to get the following "0 curvature condition" : 
\begin{equation} \label{zs+} dZ_+ -\frac{1}{2}[Z_+,Z_+] = 0 \end{equation}
which is equivalent to  
\begin{equation} \label{zs-} dZ_- -\frac{1}{2}[Z_-,Z_-] = 0 \end{equation}
As one can easily see, if $Z$ were not depending on $L$, the last equation 
would be an equation of transport on a principal bundle that should be determined. 
Unfortunately, $Z$ depends on $L$. In \cite{M1}, M. Mulase shows how this equation 
is equivalent to an equation of transport.

\subsection{The Lax equation as a gauge equation}

We have here to build up the principal bundle $P$ of the last proposition, before defining 
the (flat) connection that we shall use. 
Here, and till the end of the paper we replace
$$ \mathcal{R} = R[[t_1,t_2,...]]$$
by 
$$ \mathcal{R} = R$$
which is possible for the (trace) Fr\"olicher structure induced by the Fr\"olicher structure on $k^\N.$  
We define 
$$ G = \{ 1 + u | u \in E^{(-1)} \}. $$
This is a regular Lie group, which is endowed with the product topology \cite{ARS}. Then, since $$[L^n,L]=0$$ we have the same relations as in the last section:
$$\frac{dL}{dt_n} = [(L^n)_+,L] \Leftrightarrow  \frac{dL}{dt_n} = [-(L^n)_-,L]$$
which gives
$$ dL = [Z_-,L]$$ with $$Z_- = \sum_{n = 1}^{+\infty} -(L^n)_-dt_n . $$
Here, $Z_-$ is a connection on the principal bundle $T \times (1+E^{(-1)})$
but also a connection on the vector bundle $T \times (k\d \oplus E^{(-1)}).$
Since $$\frac{dL}{dt_n} = \frac{d(L_-)}{dt_n},$$
we get also the 0-curvature condition
$$ dZ_- - \frac{1}{2}[Z_-,Z_-] = 0$$
which gives the following, applying Theorem \ref{Ambrose-Singer}:

\begin{Proposition} 
The equation \ref{lax1} is equivalent to the following: 
\begin{equation} \label{lax2} \exists S \in C^\infty(T,G); \left\{ \begin{array}{l} 
L = S \d S^{-1}\\
dS = Z_- S \end{array} \right. \end{equation}
\end{Proposition}

$S$ is called the {\it gauge operator of $L$}, and we see that one can easily see that defining $L$ is the same as defining $S$ (up to operators that commute with $\d$).  Let us work a little more with this formulation. It gives us the following : 
\begin{eqnarray*}
dS.S^{-1} & = & \sum_{n = 1}^{\infty}dt_n\lbrace S. (\d ^n ). S^{-1} \rbrace_- \end{eqnarray*}

\subsection{Lax-type reduction of connections: the $Z_--$equation}
Since $T$ is naturally endowed with a Fr\"olicher structure componentwise,
and that $G$ is a regular Lie group \cite{ARS} in the ILH class \cite{ARS,Om}.
According to Theorem \ref{Ambrose-Singer}, if $\theta$ is a flat connection (i.e. its curvature is null), since $\pi_1(T)=0,$ thete is a global section $S$ of $P_G = T \times G$ such that $dS.S^{-1}=\theta.$
\begin{Theorem}
Let $\theta_-$ be a flat connection 1-form on $P_G.$
Then there exists a Lax operator $L$ such that $\theta_- = Z_-.$ 
\end{Theorem}

\noindent
\textbf{Proof.}
Since $T$ is naturally endowed with a Fr\"olicher structure componentwise,
and that $G$ is a regular Lie group \cite{ARS} in the ILH class \cite{ARS,Om}.
According to Theorem \ref{Ambrose-Singer}, if $\theta_-$ is a flat connection (i.e. its curvature is null), since $\pi_1(T)=0,$ there is a global section $S$ of $T \times G$ such that $dS.S^{-1}=\theta.$
Let $$ L = S \d S^{-1}.$$
Then $Z_- = dS.S^{-1} = \theta.$\qed

\subsection{The total equation, flat connections and the $Z_+$-equation}
Now, let us define
$$\tilde \A _m = \left\{ \sum_{i \in \Z, i \leq m}  q^m u_i \partial^i \right\}.$$
Here, we could understand $q$ as a formal parameter, but for the full rigor
of the sequen, we need to say that $q$ is a $k-$variable. Then, 
$$ \tilde \A =\left\{ \sum_{m = 1}^\infty a_m | a_m \in \A_m \right\}$$

and we set $${G}_\A = {G} .(1 + \A) =\left\{ g + \sum_{m = 1}^\infty a_m | g \in {G} \wedge a_m \in \A_m \right\}.$$
With $\pi : G_\A \rightarrow G$ the projection on the first component, we get that $$ (1 + \A) = Ker \pi$$
so that, we have an exact sequence of Lie groups
$$ 1 \rightarrow (1+\A) \rightarrow G_\A \rightarrow G \rightarrow 1,$$
and we have the following:
\begin{Proposition}
The Fr\"olicher Lie groups $G,$ $(1+\A)$ ans $G_\A$ are regular with regular Lie algebra.
\end{Proposition} 

\noindent
\textbf{Proof.}
The Lie group $G$ is regular from \cite{ARS}, 
the Fr\"olicher Lie group $(1+\A)$
 is regular with regular Lie algebra from 
Theorem \ref{regulardeformation} and the Fr\"olicher Lie group $G_\A$ safisfies the assumptions of Theorem \ref{exactsequence}. \qed

\vskip 12pt
We consider the 
flat bundles 
$$ P_{(1 +\A)} = T \times (1+\A) \quad \hbox{and} \quad P_{G_\A} = T \times G_\A.  $$
We consider a flat connection $\theta$ on $P_{G_\A}.$
Applying Theorem \ref{Ambrose-Singer}, we get:
\begin{Proposition}
There exits a smooth horizontal section $Q_\theta$ with respect to $\theta$ such that $Q_\theta |_{t = 0}  \in G.$
\end{Proposition}

Now, we use the projection $\pi: G_\A \rightarrow G$ and its derivative $\pi ' : E^{(-1)} \oplus \A \rightarrow E^{(-1)}$
to get that 
\begin{itemize}
\item The map $S = \pi \circ Q_\theta : T \rightarrow G$ is a smooth section of $P_G,$
\item The 1-form $\theta_- = \pi '\circ \theta$ is a flat connection on $P_G$ trivially. 
\end{itemize}
 
Then, to any (flat) connection $\theta$ we have a Lax operator $L$ such that  
$$\pi \circ \theta = \Sigma_{n = 1}^\infty (L^n)_- dt_n = Z_-.$$

We build 
$$Z_{q,+} = \sum_{n = 1}^{+\infty} q^n (L^n)_+ dt_n$$
With this, the equation of the KP hierarchy becomes
$$ {d_q (qL)} = [ Z_{q,+} , qL ] $$
with $$ d_q = \sum_{n = 1}^{+\infty} q^n dt_n = \sum_{n = 1}^{+\infty} d(q^n t_n).$$
We can see here a scaling effect due to the change of variables 
$t_n \mapsto q^nt_n,$ passing from $d$ to $d_q.$ Straight way from this, we get that 
$$ Z_{q,+} = S.\omega_q .S^{-1} + d_q S . S^{-1} $$
with $\omega_q = \sum_{n = 1}^{+\infty} q^n \d^n$
and that \begin{equation} \label{zsq+} d_q Z_{q,+} - \frac{1}{2}[Z_{q,+},Z_{q,+}] = 0 \end{equation}
which is a $0-$curvature condition. 
 
In other words, with the section $S$ obtained in Equation (\ref{lax2}), we define the smooth section $$\tilde{S}(t_1,t_2,...) = S(qt_1, q^2t_2,...)$$
and the equation (\ref{lax1}) is equivalent to 
$$ d \tilde L = [Z_{q,+}, \tilde L] = [Z_{q,-}, \tilde L]$$
with $$ Z_{q,-} = d_qS.S^{-1} = d \tilde{S} . \tilde{S}^{-1},$$
 $$ \tilde{L} = q L(qt_1, q^2t_2,...) = \tilde{S} . q\d . \tilde{S}^{-1} $$
and 
$$   Z_{q,+} = \sum_{n = 1}^{+\infty} (\tilde{L}^n)_+ dt_n = \tilde{S}. \left( \sum_{n = 1}^{+\infty} q^n \d^n\right) .\tilde{S}^{-1} + d \tilde{S} . \tilde{S}^{-1}.$$
Equation (\ref{zs+}) is then exactly the same as (\ref{zsq+}) and through the change of variables announced, we get $$dZ_{q,+} -\frac{1}{2}[Z_{q,+},Z_{q,+}]=0.$$
Thus, $Z_{q,+}$ is a $\A-$valued connection with $0-$curvature. Then, applying Theorem \ref{Ambrose-Singer}, we get
\begin{Theorem} \label{z+}
There exists an unique global smooth section $Y_q$ of $P_{1+\A}$ such that $dY_q.Y_q^{-1} = Z_{q,+}$ and $Y_q(0)=1.$  
\end{Theorem}  
Finally, we finish with the analog of \cite{M1}, Theorem 1:
\begin{Theorem}
Let $U_q = \tilde{S}^{-1} Y_q.$
Then the Lax equation \ref{lax1} of the KP hierarchy is equivalent to 
$$ dU_q = \omega_q U_q.$$
\end{Theorem}

\noindent
\textbf{Proof.}
We have on one hand: 
\begin{eqnarray*}
\tilde{S} \omega_q \tilde{S}^{-1} & = & \sum_{n = 1}^{+\infty} \tilde{L}^n dt_n\\
& = & \sum_{n = 1}^{+\infty} (\tilde{L}^n)_+ dt_n +\sum_{n = 1}^{+\infty} (\tilde{L}^n)_- dt_n\\
& = & Z_{q,+} - d\tilde{S}.\tilde{S}^{-1} \\
& = & dY_q.Y_q^{-1} - d\tilde{S}.\tilde{S}^{-1}
\end{eqnarray*}
and on the other hand:
\begin{eqnarray*}
dU_q & = & d\left(\tilde{S}^{-1}.Y_q\right)\\
& = & - \tilde{S}^{-1} d\tilde{S}. \tilde{S}^{-1} .Y_q + \tilde{S}^{-1} dY_q
\end{eqnarray*}
Thus 
\begin{eqnarray*} \omega_q U_q & = & \tilde{S}^{-1} . \left(dY_q.Y_q^{-1} - d\tilde{S}.\tilde{S}^{-1}\right) Y_q \\
& = &\tilde{S}^{-1} dY_q - \tilde{S}^{-1} d\tilde{S}. \tilde{S}^{-1} .Y_q \\
& = & dU_q \qed \end{eqnarray*}
\subsection{The Fr\"olicher Lie group of formal sections: last remarks}

Now, in the approach by formal sections, smooth sections of the bundles $P_G,$
$P_{G_\A}$ and $P_{(1+\A)}$ are understood as formal power series (formally analytic functions) in the variables $t_1,t_2,...$ or can be viewed as smooth germs of mappings at $t=0.$  This last approach by germs appears to us more adapted with the results of last section, espacially from the viewpoint of Mulase's paper \cite{M1} where the (formal) power series in $(t_1,t_2...)$ are graded with respect to the orders of the monomials $t_1^{k_1}t_2^{k_2}...$ (finite products) with the valuations $val(t_n) = n.$
From the Fr\"olicher Lie group $C^\infty(T,G)$ of (true) smooth sections
of the bundle $P_G,$ we get the Fr\"olicher Lie group analogous to the infinite jet space at $t=0$ usually noted as ``$J^\infty_0(T,G)$''.  
We can rewrite this relation as $val(t_n) = val(q^n)$ at the 
level of germs of sections and recover the finite sums 
in our graded algebra $\A.$ 
This valuation is then no longer algebraic, but is the key tool to 
get regular Lie groups (with exponential maps) at the level of sections.   
Finally, we remark that Mulase's Birkoff decomposition \cite{M1}
on germs of sections has an analog on Lie groups of sections on the 
short exact sequence $$ 1 \rightarrow (1+\A) \rightarrow G_\A \rightarrow G \rightarrow 1,$$
and that our analysis of the KP equations shows that we get a section on
the bundle $T \times (1+\A)_+$ with the map $Y_q$ obtained in Theorem \ref{z+} .


\begin{thebibliography}{99}


\bibitem{ARS} M. Adams, T.Ratiu, R. Schmid {\it A Lie group structure for pseudodifferential operators}
Math. Ann. \textbf{273} no4 (1986), 529-551

\bibitem{CN} Cherenack, P.; Ntumba, P.; Spaces with differentiable structure an application to cosmology 
\textit{Demonstratio Math.} \textbf{34} no 1 (2001), 161-180 


\bibitem{Don} Donato, P.; \textit{Rev\^etements de groupes diff\'erentiels} Th\`ese de doctorat d'\'etat, 
Universit\'e de Provence, Marseille (1984) 


\bibitem{FK} Fr\"olicher, A; Kriegl, A; {\it Linear spaces and differentiation theory} Wiley series in Pure and Applied Mathematics, Wiley Interscience (1988)


\bibitem{Gil} Gilkey, P;
{\it Invariance theory, the heat equation and the Atiyah-Singer index theorem}
Publish or Perish (1984)

\bibitem{Glo} Gl\"ockner, H;
Algebras whose groups of the units are Lie groups
\textit{Studia Math. } \textbf{153}, no2 (2002), 147-177

\bibitem{Ig} Iglesias, P.; Connexions et diff\'eologie 
\textit{Aspects dynamiques et topologiques des groupes infinis de transformation de la m\'ecanique}
Travaux en cours \textbf{25}, Hermann (1987), 61-78

\bibitem{Igdiff} Iglesias, P.; \textit{Diffeology} http://math.huji.ac.il/~piz/documents/Diffeology.pdf

\bibitem{KMS} Kolar, I.; Michor, P.W.; Slovak, J.; \textit{Natural operations in differential geometry}; Springer (1993)

\bibitem{KM} Kriegl, A.; Michor, P.W.; \textit{The convenient setting for global analysis} 
Math. surveys and monographs \textbf{53}, American Mathematical society, Providence, USA. (2000)



\bibitem{L} Lax, P. D.; {\it Integrals of non linear equations of evolution and solitary waves} Comm. Pure Appl. Math. {\bf 21} p. 467-490 (1968)

\bibitem{Les} Leslie, J.; On a Diffeological Group Realization of certain 
Generalized symmetrizable Kac-Moody Lie Algebras
\textit{J. Lie Theory} \textbf{13} (2003), 427-442

\bibitem{Li} Lichnerowicz, A.  {\it Th\'eorie globale des connexions et des groupes d'holonomie}  ed. Cremonese, Roma (1956)

\bibitem{Ma} Magnot, J-P.; Diff\'eologie du fibr\'e d'Holonomie en dimension infinie, \textit{ Math. Rep. Can. Roy. Math. Soc.} \textbf{28} no4 (2006) 

\bibitem{Ma2} Magnot, J-P.; Structure groups and holonomy in infinite dimensions, \textit{Bull. Sci. Math.}
\textbf{128} (2004), 513-529


\bibitem{M1} Mulase, M.; {Complete integrability of the Kadomtsev-Petviashvili equation} \textit{Adv. Math.} {\bf 54} (1984) 57-66 

\bibitem{M2} Mulase, M.; { Solvability of the super KP equation and a generalization of the Birkhoff decomposition} \textit{Invent. Math.} {\bf 92} (1988), 1-46 

\bibitem{M3} Mulase, M.; Geometry of soliton equations; MSRI preprint (1983)



\bibitem{Om} Omori, H.; \textit{Infinite dimensional Lie groups} AMS translations of mathematical monographs \textbf{158} (1997)

\bibitem{Sou} Souriau, J.M.; Un algorithme g\'en\'erateur de structures quantiques; 
\textit{Ast\'erisque}, Hors S\'erie, (1985) 341-399 

\bibitem{Rob} Robart, T.; Sur l'int\'egrabilit\'e des sous-alg\`ebres de Lie en dimension infinie; \textit{Can. J. Math.} \textbf{49} (4) (1997), 820-839

\bibitem{Ue} Ueno, K.; Analytic and algebraic aspects of the Kadomtsev-Petviashvili hierarchy from the viewpoint of the universal Grassmann manifold; \textit{ Infinite Dimensiona groups with applications, MSRI publications} \textbf{4}, Springer (1985), 335-353

\bibitem{Wa} Watts, J.; \textit{Diffeologies, differentiable spaces and symplectic geometry} PhD thesis arXiv:1208.3634v1

\end{thebibliography}
\end{document}